\documentclass[11pt,letterpaper]{article}
\usepackage[margin=1in,headheight=14pt]{geometry}
\usepackage[T1]{fontenc}
\usepackage{lmodern}
\usepackage{amsmath,amssymb,amsthm,booktabs,array,enumitem}
\usepackage{fancyhdr}
\usepackage[colorlinks=true,linkcolor=blue!45!black,citecolor=blue!45!black,urlcolor=blue!45!black]{hyperref}
\usepackage{xcolor}
\hypersetup{pdftitle={A Lower Bound for the Chromatic Number of Ball-Packing Graphs},pdfauthor={Dan Ismailescu; Sunwoo Lee; Terry Lee}}
\setlist[itemize]{topsep=3pt,itemsep=2pt}
\setlist[enumerate]{topsep=3pt,itemsep=3pt}
\allowdisplaybreaks[1]
\numberwithin{equation}{section}
\newtheorem{theorem}{Theorem}[section]
\newtheorem{proposition}[theorem]{Proposition}
\newtheorem{lemma}[theorem]{Lemma}

\theoremstyle{definition}
\newtheorem{definition}[theorem]{Definition}
\newtheorem{remark}[theorem]{Remark}
\newcommand{\F}{\mathbb F}
\newcommand{\R}{\mathbb R}
\newcommand{\one}{\mathbf1}
\newcommand{\im}{\operatorname{im}}
\newcommand{\rank}{\operatorname{rank}}

\newcommand{\tr}{\operatorname{tr}}
\newcommand{\chiball}{\chi_{\mathrm{ball}}}
\newcommand{\chicong}{\chi_{\mathrm{cong}}}
\title{A Lower Bound for the Chromatic Number\\of Ball-Packing Graphs}
\author{%
Dan Ismailescu\textsuperscript{1}\thanks{Corresponding author: \href{mailto:dan.p.ismailescu@hofstra.edu}{\texttt{dan.p.ismailescu@hofstra.edu}}.}\quad
Sunwoo Lee\textsuperscript{2}\quad
Terry Lee\textsuperscript{3}\\[0.6em]
\small\textsuperscript{1}Mathematics Department, Hofstra University, Hempstead, NY\\
\small\textsuperscript{2}Taft School, Watertown, CT\\
\small\textsuperscript{3}Thomas Jefferson School, St.\ Louis, MO
}
\date{}
\begin{document}
\maketitle
\vspace{-2em}
\begin{abstract}
For every odd prime power $q$, we construct a finite packing of congruent balls in dimension $d=q^2-q+1$ whose contact graph requires at least $2d$ colors. Consequently, the ball-packing chromatic number satisfies $\chiball(q^2-q+1)\ge2(q^2-q+1)$, and the same lower bound holds for congruent-ball packings. The construction starts from the complements of Hubaut's $C_{20}$ graphs and uses their Taylor extensions. We give an explicit finite-field adjacency rule for the $q^3$-vertex graph $H_q$, prove its strong regularity, determine its spectrum and projection matrix, establish $\alpha(H_q)=q$, and construct a two-distance embedding in $\R^{q(q-1)}$. Its Taylor extension $T_q$ has $2(q^3+1)$ vertices. We prove its distance-regularity, compute its intersection array, construct a three-distance embedding in $\R^{d}$, and prove $\alpha(T_q)=q+1$. The chromatic lower bound follows by realizing $T_q$ as the contact graph of equal-radius balls and bounding the size of every color class. All coordinates, distances, and independence bounds are derived explicitly. The case $q=3$ uses the Schl\"afli and Gosset graphs and gives a packing of 56 congruent balls in $\R^7$ requiring at least 14 colors.
\end{abstract}

\section{Introduction and conventions}
The graphs considered here belong to the classical connection between strongly regular graphs, regular two-graphs, and antipodal distance-regular graphs. Hubaut's survey~\cite{Hubaut1975} is an early systematic account of strongly regular graphs. Taylor's work~\cite{Taylor1977} develops regular two-graphs and their relation to the constructions underlying the present family. Modern treatments appear in Brouwer and Van Maldeghem~\cite{BVM2022} and, for distance-regular graphs, in Brouwer, Cohen, and Neumaier~\cite{BCN1989}.

The exposition is intended to be accessible to a nonspecialist in graph theory while remaining comprehensive. We give a single explicit route from a finite-field adjacency test to real Euclidean coordinates, with the intervening proofs included. We distinguish carefully between finite-field coordinates, which define the graph, and real coordinates, which realize its vertices as points.

\subsection{The ball-packing consequence}
A ball packing in $\R^m$ is a collection of closed Euclidean balls of positive radii with pairwise disjoint interiors. Its contact graph has one vertex for each ball and an edge precisely when two balls are tangent. Define
\[
\chiball(m)=\sup_{\mathcal P}\chi(G(\mathcal P)),
\qquad
\chicong(m)=\sup_{\mathcal P\text{ congruent}}\chi(G(\mathcal P)),
\]
where the first supremum is over finite ball packings in $\R^m$, and the second is over finite packings in which all radii are equal. A common dilation turns a congruent-ball packing into a packing of unit-radius balls without changing its contact graph. In particular, $\chicong(m)\le\chiball(m)$.

The packing consequence of our construction is the following theorem. It is proved in Section~\ref{sec:packingproof} from the explicit distances and the exact independence number of the Taylor extension.
\begin{theorem}[Ball-packing chromatic lower bound]\label{thm:packing}
For every power $q$ of an odd prime, put $d=q^2-q+1$. There is an explicit finite packing of $2(q^3+1)$ congruent balls in $\R^{d}$ whose contact graph requires at least $2d$ colors. Consequently,
\begin{equation}\label{eq:packing-intro}
\boxed{\chiball(q^2-q+1)\ \ge\ \chicong(q^2-q+1)\ \ge\ 2(q^2-q+1).}
\end{equation}
\end{theorem}
For example, the dimensions $7,21,43,73$ obtained from $q=3,5,7,9$ admit congruent-ball packings requiring at least $14,42,86,146$ colors, respectively.

\subsection{Related work on coloring ball packings}
Maehara~\cite{Maehara2007} studied chromatic numbers of solid-ball configurations in three dimensions, including both arbitrary balls and unit balls. Bagchi and Datta~\cite{BagchiDatta2013} formulated higher-dimensional ball-packing coloring as an analogue of the map-coloring problem. These papers concern contact graphs, in which interiors are disjoint, rather than general intersection graphs. The distinction between congruent balls and balls of varying radii is also relevant; Theorem~\ref{thm:packing} applies to the congruent subclass.

Chen~\cite{Chen2017} used Euclidean representations of strongly regular graphs to construct ball packings with large chromatic number. In particular, for odd prime powers $q$, he considered the complements of Hubaut's $C_{20}$ family, realized in dimension $q^2-q$ with chromatic number at least $q^2$. He also obtained congruent-ball packings in dimensions $q^3-q^2+q$, for all prime powers $q$, whose contact graphs have chromatic number $q^3+1$. These are different dimension sequences from the one in Theorem~\ref{thm:packing}. Chen's approach was motivated by the use of strongly regular two-distance sets in Bondarenko's work on Borsuk's conjecture~\cite{Bondarenko2014}.

The graph-theoretic and spectral ingredients used here are classical. We present their explicit combination with the Taylor extension and isolate the resulting inequality for $\chiball$ and $\chicong$. The proof uses the contact graph $T_q$, its $2(q^3+1)$ vertices, its independence number $q+1$, and its realization in dimension $q^2-q+1$. No optimality claim for these packing chromatic bounds is needed.

\subsection{Graph conventions and main structural theorem}
There are two complementary adjacency conventions in the literature. We use $H_q$ for the graph of degree $(q+1)(q^2-1)/2$. Its complement is the family commonly denoted $C_{20}(q)$; it has degree $(q-1)(q^2+1)/2$. In particular, $H_3$ is the degree-$16$ Schl\"afli graph, rather than its degree-$10$ complement. The spectral formulas for this convention also appear in Mohar~\cite{Mohar2009}. The distinction matters when taking a Taylor extension: the convention used here gives the Gosset graph when $q=3$. Antipodal covers of complete graphs and their strongly regular neighborhoods are also treated by Godsil~\cite{Godsil1992}.

All graphs are finite, simple, and undirected. An independent set contains no adjacent pair; its largest size is $\alpha(G)$. The chromatic number $\chi(G)$ is the least number of colors in a proper vertex coloring. A graph is strongly regular with parameters $(n,k,\lambda,\mu)$ if it has $n$ vertices, every degree is $k$, adjacent vertices have $\lambda$ common neighbors, and distinct nonadjacent vertices have $\mu$ common neighbors.

For the remainder of the paper, set
\begin{equation}\label{eq:parameters}
\begin{aligned}
n&=q^3,&
k&=\frac{(q+1)(q^2-1)}2,\\
\lambda&=\frac{q^3+3q^2-3q-5}{4},&
\mu&=\frac{(q+1)(q^2-1)}4,\\
f&=q(q-1),&
B&=\frac{(q-1)(q^2+1)}2=n-k-1.
\end{aligned}
\end{equation}
Here $B$ is a number, not a matrix. The matrices $I,J$ always have the size dictated by context; $J$ has every entry equal to one.

\begin{theorem}\label{thm:main}
For every power $q$ of an odd prime, the construction below gives the following graphs.
\begin{enumerate}[label=\textup{(\alph*)}]
\item $H_q$ is strongly regular with parameters $(n,k,\lambda,\mu)$ in~\eqref{eq:parameters}. Its spectrum is
\[
k^1,\qquad
\left(\frac{q^2-1}{2}\right)^{q(q-1)},\qquad
\left(-\frac{q+1}{2}\right)^{(q-1)(q^2+1)}.
\]
\item If $A$ is its adjacency matrix, the orthogonal projection onto the $f$-dimensional nonconstant eigenspace is
\[
E=\frac{2}{q(q+1)}A+\frac1qI-\frac1{q^2}J.
\]
Moreover, $\alpha(H_q)=q$, and $H_q$ has a two-distance embedding in $\R^f$ with squared distances
\[
\frac{2(q-1)}{q(q+1)}\quad\text{on edges},
\qquad \frac2q\quad\text{on distinct nonedges}.
\]
\item The Taylor extension $T_q$ has $2(n+1)$ vertices, degree $n$, and intersection array
\[
\{n,B,1;\ 1,B,n\}.
\]
Its independence number is $q+1$. It has a three-distance embedding in $\R^{f+1}$, with squared distances
\[
\frac{2(q-1)}{q(q+1)},\qquad \frac2q,\qquad\frac4{q+1}
\]
for pairs at graph distances $1,2,3$, respectively.
\end{enumerate}
\end{theorem}
The superscripts in the displayed spectrum denote multiplicities. We prove every assertion of the theorem below.

\section{Finite-field arithmetic and squares}
\subsection{The field of \texorpdfstring{$q$}{q} elements}
Write $q=p^m$, with $p$ an odd prime. We work in a field $F=\F_q$. If $m=1$, this is arithmetic modulo $p$. For computations with $m>1$, choose an irreducible polynomial $g(t)$ of degree $m$ over $\F_p$ and form
\[
F=\F_p[t]/(g(t)).
\]
Every element is represented by a polynomial of degree less than $m$. Add coefficients modulo $p$; multiply polynomials and reduce the result modulo $g(t)$. This is field arithmetic, rather than arithmetic modulo the composite integer $p^m$.

For completeness, existence of a field with $q$ elements follows from the roots of $X^q-X$ in a splitting field over $\F_p$. Its derivative is $-1$, so it has $q$ distinct roots. In characteristic $p$, the identity $(a+b)^q=a^q+b^q$ shows that its roots are closed under addition and subtraction; they are also closed under multiplication and nonzero inversion. They therefore form a field of size $q$. A finite splitting field can be obtained by successively adjoining roots of irreducible factors. The polynomial-quotient description above is the standard computational model.

\subsection{Squares and the quadratic character}
We write $F^\times=F\setminus\{0\}$ for the set of nonzero elements of $F$. Thus $u\in F^\times$ means $u\in F$ and $u\ne0$; in particular, division by $u$ is allowed. Under multiplication, $F^\times$ is called the multiplicative group of $F$. The same notation is used for other fields below.

A nonzero element of $F$ is a square if it equals $z^2$ for some $z\ne0$. Since $z^2=w^2$ implies $(z-w)(z+w)=0$, every nonzero square has precisely two square roots, $z$ and $-z$. They are different because the characteristic is odd. Thus there are $(q-1)/2$ nonzero squares and the same number of nonsquares.

Define the quadratic character by
\[
\chi(a)=\begin{cases}
0,&a=0,\\
1,&a\text{ is a nonzero square},\\
-1,&a\text{ is a nonsquare}.
\end{cases}
\]
The nonzero squares form a subgroup of index two in $F^\times$. Its two cosets give
\begin{equation}\label{eq:character-product}
\chi(ab)=\chi(a)\chi(b).
\end{equation}
This includes zero factors. Fix a nonsquare $\delta{}\in F$ for the construction.

\begin{lemma}\label{lem:charsum}
For every nonzero $r\in F$,
\[
\sum_{s\in F}\chi(s^2-r)=-1.
\]
\end{lemma}
\begin{proof}
Count pairs $(s,t)\in F^2$ satisfying $t^2=s^2-r$. For a fixed $s$, the number of possible $t$ is $1+\chi(s^2-r)$: it is one, two, or zero according as $s^2-r$ is zero, a nonzero square, or a nonsquare.

On the other hand, the equation is $(s-t)(s+t)=r$. Choose $u=s-t$ arbitrarily in $F^\times$; then $s+t=r/u$, so
\[
s=\frac{u+r/u}{2},\qquad t=\frac{r/u-u}{2}.
\]
Each of the $q-1$ choices gives exactly one pair. Hence
\[
q+\sum_s\chi(s^2-r)=q-1,
\]
which proves the identity.
\end{proof}

\section{The graph on triples}
\begin{definition}\label{def:H}
The vertex set of $H_q$ is $F^3$. For distinct vertices $P=(a,b,c)$ and $Q=(A,B,C)$, put
\begin{equation}\label{eq:DS}
D=(a-A)^2-\delta{}(b-B)^2,\qquad S=c-C+bA-aB.
\end{equation}
Declare
\begin{equation}\label{eq:adjacency}
P\sim Q\quad\Longleftrightarrow\quad
4\delta{}S^2-D^2\text{ is a nonzero square in }F.
\end{equation}
\end{definition}
Every operation in~\eqref{eq:DS}--\eqref{eq:adjacency} takes place in $F$. There are $q^3$ vertices. Exchanging $P,Q$ leaves $D$ unchanged and sends $S$ to $-S$, so adjacency is symmetric.

\begin{lemma}\label{lem:anisotropic}
For $u,v\in F$, the equality $u^2-\delta{}v^2=0$ implies $u=v=0$.
\end{lemma}
\begin{proof}
If $v\ne0$, the equality gives $\delta{}=(u/v)^2$, contradicting the choice of $\delta{}$. If $v=0$, it gives $u^2=0$, hence $u=0$.
\end{proof}

\begin{proposition}\label{prop:degree}
The graph $H_q$ has degree $k=(q+1)(q^2-1)/2$. For fixed $a,b$, the set
\[
I_{a,b}=\{(a,b,c):c\in F\}
\]
is independent. These $q^2$ sets partition the vertices.
\end{proposition}
\begin{proof}
Fix $P=(a,b,c)$. If $(A,B)=(a,b)$ and $Q\ne P$, then $D=0$ and $S=c-C\ne0$. The edge-test expression is $4\delta{}(c-C)^2$, a nonsquare. This proves independence of each $I_{a,b}$ and excludes neighbors in the same part.

For any of the other $q^2-1$ choices of $(A,B)$, Lemma~\ref{lem:anisotropic} gives $D\ne0$. As $C$ ranges over $F$, the quantity $S=c-C+bA-aB$ ranges over $F$ exactly once. For such fixed $D$,
\begin{align*}
\sum_{s\in F}\chi(4\delta{}s^2-D^2)
&=\chi(4\delta{})\sum_{s\in F}\chi\left(s^2-\frac{D^2}{4\delta{}}\right)\\
&=(-1)(-1)=1
\end{align*}
by Lemma~\ref{lem:charsum}. Also $4\delta{}s^2-D^2$ is never zero: if it were, $s\ne0$ and $\delta{}=(D/(2s))^2$ would be a square. Among these $q$ character values, therefore, exactly $(q+1)/2$ are $+1$. Each of the $q^2-1$ choices of $(A,B)$ contributes this many neighbors. Multiplying proves the degree formula.
\end{proof}

\begin{remark}[The choice of nonsquare]
The isomorphism type does not depend on which nonsquare $\delta{}$ is chosen. If $\delta{}'=\delta{}t^2$ for $t\in F^\times$, the bijection
\[
(a,b,c)\longmapsto(a,b/t,c/t)
\]
preserves the edge-test expression: $D$ stays the same, $S$ is divided by $t$, and $\delta{}'S'^2=\delta{}S^2$. Any two nonsquares have a square ratio because they belong to the same coset of the subgroup of squares.
\end{remark}

\section{An elementary proof of strong regularity}
The degree calculation alone does not prove strong regularity: common-neighbor counts must also be independent of the chosen pair. We derive those counts from an explicit matrix identity. The auxiliary field and transformations in this section are used to prove that identity; the graph itself is already completely specified by Definition~\ref{def:H}.

\subsection{A quadratic extension and its norm}
Adjoin a symbol $\theta$ satisfying $\theta^2=\delta{}$. Because $\delta{}$ is a nonsquare, $X^2-\delta{}$ is irreducible over $F$, and
\[
L=F[\theta]=\{a+b\theta:a,b\in F\}
\]
is a field. Define conjugation and norm by
\[
\overline{a+b\theta}=a-b\theta,
\qquad
N(a+b\theta)=(a+b\theta)(a-b\theta)=a^2-\delta{}b^2.
\]
Conjugation respects addition and multiplication and applying it twice does nothing. Consequently
\[
N(zw)=N(z)N(w),\qquad N(\bar z)=N(z).
\]
Lemma~\ref{lem:anisotropic} shows that $N(z)=0$ only when $z=0$. These identities can also verify division explicitly: for $z\ne0$, its inverse is $\bar z/N(z)$.

Put
\begin{equation}\label{eq:eta}
\eta(z)=\chi(N(z)),\qquad \varepsilon=\chi(-1).
\end{equation}
For nonzero $z$, $\eta(z)$ is $+1$ or $-1$, and
\begin{equation}\label{eq:eta-properties}
\eta(zw)=\eta(z)\eta(w),\qquad
\eta(\bar z)=\eta(z).
\end{equation}
We use only these identities, not any classification of squares in $L$.

\subsection{Representatives for the vertices and an extra point}
Represent $P=(a,b,c)$ by
\begin{equation}\label{eq:representatives}
r_P=\left(-\frac{a^2-\delta{}b^2}{2}+c\theta,\ a+b\theta,\ 1\right)\in L^3.
\end{equation}
Introduce one more representative $r_\infty=(1,0,0)$. This extra index is an aid to the proof and is not a vertex of $H_q$.

For $w=(x,y,z)$ and $w'=(X,Y,Z)$ in $L^3$, define
\begin{equation}\label{eq:hermitian}
h(w,w')=x\bar Z+y\bar Y+z\bar X.
\end{equation}
This expression is linear in its first argument and conjugate-linear in the second, and $h(w',w)=\overline{h(w,w')}$. Such an expression is called a Hermitian form, but we will use only its displayed formula.

Every chosen representative satisfies $h(r_i,r_i)=0$. In fact, they represent all nonzero solutions of this equation up to nonzero scalar multiplication. To verify this, let $w=(x,y,z)\ne0$ satisfy
\[
x\bar z+N(y)+z\bar x=0.
\]
If $z=0$, then $N(y)=0$, so $y=0$ and $w$ is a nonzero scalar multiple of $r_\infty$. If $z\ne0$, divide all entries by $z$. The normalized vector has the form $(x',y',1)$ and satisfies
\[
x'+\bar x'+N(y')=0.
\]
Writing $y'=a+b\theta$ and $x'=u+c\theta$, the equation becomes $2u+a^2-\delta{}b^2=0$. Thus $u=-(a^2-\delta{}b^2)/2$, giving exactly~\eqref{eq:representatives}. The representative is unique because its last entry is one.

For two finite vertices, expansion gives
\begin{equation}\label{eq:hDS}
h(r_P,r_Q)=-\frac D2+S\theta.
\end{equation}
Indeed its real part is
\[
-\frac{a^2-\delta{}b^2}{2}-\frac{A^2-\delta{}B^2}{2}+aA-\delta{}bB
=-\frac{(a-A)^2-\delta{}(b-B)^2}{2},
\]
and its $\theta$-coefficient is $c-C+bA-aB$. Therefore
\begin{equation}\label{eq:sign-equivalence}
\varepsilon\eta(h(r_P,r_Q))
=\chi\left(\delta{}S^2-\frac{D^2}{4}\right)
=\chi(4\delta{}S^2-D^2).
\end{equation}
The last equality uses $\chi(4)=1$.

For distinct representatives, $h(r_i,r_j)\ne0$. This is immediate if one index is $\infty$, since $h(r_\infty,r_P)=1$. For finite indices, zero in~\eqref{eq:hDS} would give $D=S=0$. Lemma~\ref{lem:anisotropic} would then give $a=A,b=B$, followed by $c=C$, contradicting distinctness.

\subsection{Explicit transformations preserving the form}
Consider
\begin{equation}\label{eq:W}
W(x,y,z)=(z,y,x),
\end{equation}
and, for $u,t\in L$ with $t+\bar t+N(u)=0$, consider
\begin{equation}\label{eq:T}
T_{u,t}(x,y,z)=(x-\bar u\,y+tz,\ y+uz,\ z).
\end{equation}
Both transformations are invertible and preserve $h$. The assertion for $W$ follows by swapping the first and third terms of~\eqref{eq:hermitian}. For $T_{u,t}$, substitute its coordinates to obtain
\begin{align*}
h(T_{u,t}w,T_{u,t}w')
={}&h(w,w')-\bar u\,y\bar Z-u z\bar Y
+t z\bar Z+\bar t z\bar Z\\
&\quad+\bar u\,y\bar Z+u z\bar Y+N(u)z\bar Z\\
={}&h(w,w')+(t+\bar t+N(u))z\bar Z=h(w,w').
\end{align*}
The map $T_{u,t}$ is invertible because $z$ is recovered first, then $y$, and then $x$.

Let $r_0=(0,0,1)$. Every finite representative $(t,u,1)$ satisfies $t+\bar t+N(u)=0$, and
\[
T_{u,t}r_0=(t,u,1),\qquad Wr_0=r_\infty.
\]
Thus $WT_{u,t}^{-1}$ moves any prescribed finite representative to $r_\infty$. Every form-preserving invertible transformation permutes the nonzero solutions of $h(w,w)=0$ up to scalar multiplication. By the normalization above, it therefore permutes our $n+1$ indices, possibly multiplying individual representatives by nonzero scalars. No general result about transformation groups is needed here.

\subsection{The sign-matrix identity}
Define a real symmetric matrix $\mathcal C$ on the $n+1$ indices by
\begin{equation}\label{eq:Cdefinition}
\mathcal C_{ii}=0,\qquad
\mathcal C_{ij}=\varepsilon\eta(h(r_i,r_j))\quad(i\ne j).
\end{equation}
All its off-diagonal entries are $\pm1$. Symmetry follows from conjugation invariance in~\eqref{eq:eta-properties}.

Let $A$ be the adjacency matrix of $H_q$, and put
\[
K=2A-J+I,\qquad \tau=q^2-q.
\]
The matrix $K$ has zero diagonal, $+1$ on edges, and $-1$ on distinct nonedges. By~\eqref{eq:sign-equivalence},
\begin{equation}\label{eq:Cblock}
\mathcal C=\begin{pmatrix}
0&\varepsilon\one^T\\
\varepsilon\one&K
\end{pmatrix}.
\end{equation}
The degree calculation gives
\begin{equation}\label{eq:Krows}
K\one=(2k-(n-1))\one=\tau\one.
\end{equation}
Hence for each finite index $P$,
\begin{equation}\label{eq:Cfirstrow}
(\mathcal C^2)_{\infty P}
=\varepsilon\sum_QK_{QP}
=\varepsilon\tau
=\tau\mathcal C_{\infty P}.
\end{equation}

We next transfer this equality to every off-diagonal entry. Let a form-preserving invertible transformation $M$ send
\[
Mr_i=\alpha_i r_{\pi(i)},\qquad \alpha_i\in L^\times,
\]
where $\pi$ is the induced permutation. Since
\[
h(r_i,r_j)=h(Mr_i,Mr_j)
=\alpha_i\bar\alpha_jh(r_{\pi(i)},r_{\pi(j)}),
\]
the signs $\delta_i=\eta(\alpha_i)\in\{1,-1\}$ satisfy
\begin{equation}\label{eq:switch}
\mathcal C_{\pi(i),\pi(j)}=\delta_i\delta_j\mathcal C_{ij}.
\end{equation}
The same transformation rule holds for the square of the matrix. Indeed, reindexing the sum by $\pi$ and using $\delta_\ell^2=1$ gives
\begin{align*}
(\mathcal C^2)_{\pi(i),\pi(j)}
&=\sum_\ell\mathcal C_{\pi(i),\pi(\ell)}\mathcal C_{\pi(\ell),\pi(j)}\\
&=\sum_\ell\delta_i\delta_\ell\mathcal C_{i\ell}
\delta_\ell\delta_j\mathcal C_{\ell j}
=\delta_i\delta_j(\mathcal C^2)_{ij}.
\end{align*}
For any $i\ne j$, choose $M$ moving $i$ to $\infty$. Equation~\eqref{eq:Cfirstrow} applies to the transformed pair. The two transformation rules then give $(\mathcal C^2)_{ij}=\tau\mathcal C_{ij}$ for the original pair.

Each diagonal entry of $\mathcal C^2$ equals $n$, because a row has $n$ off-diagonal entries with square one. We have proved
\begin{equation}\label{eq:Cidentity}
\boxed{\mathcal C^2=nI+\tau\mathcal C.}
\end{equation}
Taking the lower-right block in~\eqref{eq:Cblock} yields
\begin{equation}\label{eq:Kidentity}
\boxed{K^2=nI+\tau K-J.}
\end{equation}

\subsection{The common-neighbor counts}
Substitute $K=2A-J+I$ into~\eqref{eq:Kidentity}. Since $AJ=JA=kJ$ and $J^2=nJ$,
\[
K^2=4A^2+4A+I+(n-4k-2)J.
\]
The right side of~\eqref{eq:Kidentity} is
\[
(n+\tau)I+2\tau A-(\tau+1)J.
\]
Using $2k=n+\tau-1$ and equating the two expressions gives
\begin{equation}\label{eq:Aidentity}
\boxed{A^2=\frac k2 I+\frac{\tau-2}{2}A+\frac k2J.}
\end{equation}
For distinct $x,y$, the entry $(A^2)_{xy}=\sum_zA_{xz}A_{zy}$ counts their common neighbors. It equals $k/2$ on a nonedge and $k/2+(\tau-2)/2$ on an edge. These values are precisely
\[
\mu=\frac{(q+1)(q^2-1)}4,
\qquad
\lambda=\frac{q^3+3q^2-3q-5}{4}.
\]
Together with Proposition~\ref{prop:degree}, this proves strong regularity and $k=2\mu$.

\section{Eigenvalues and multiplicities}
The adjacency matrix $A$ is real symmetric. Its row sums are $k$, so $A\one=k\one$. Let
\[
\one^\perp=\{z\in\R^n:\one^Tz=0\}.
\]
This subspace is preserved by $A$, since $\one^TAz=(A\one)^Tz=k\one^Tz=0$. Also $Jz=0$ on this subspace. Thus~\eqref{eq:Aidentity} becomes
\[
\left(A-\frac{q^2-1}{2}I\right)
\left(A+\frac{q+1}{2}I\right)z=0
\qquad(z\in\one^\perp).
\]
Put
\begin{equation}\label{eq:rs}
r=\frac{q^2-1}{2},\qquad s=-\frac{q+1}{2}.
\end{equation}
The subspace $\one^\perp$ is the direct sum of these two eigenspaces. One can see this without using a diagonalization theorem: for $z\in\one^\perp$, define
\[
z_r=\frac{(A-sI)z}{r-s},\qquad
z_s=\frac{(rI-A)z}{r-s}.
\]
Then $z=z_r+z_s$, and the factored identity gives $Az_r=rz_r$, $Az_s=sz_s$. The two eigenspaces are perpendicular. Indeed, symmetry gives
\[
r(u\cdot v)=(Au)\cdot v=u\cdot(Av)=s(u\cdot v),
\]
so $u\cdot v=0$ when $r\ne s$.

Because $k$ differs from $r,s$, the $k$-eigenspace is exactly the one-dimensional constant space. If $f,g$ are the dimensions of the $r$- and $s$-eigenspaces, then
\[
f+g=n-1,\qquad k+rf+sg=\tr A=0.
\]
Here the trace is zero because $A$ has zero diagonal. The trace equals the sum of eigenvalues in an eigenvector basis; its invariance under a change of basis follows from $\tr(UV)=\tr(VU)$.

Solving gives
\begin{equation}\label{eq:multiplicities}
f=\frac{-k-s(n-1)}{r-s}=q(q-1),
\qquad g=n-1-f=(q-1)(q^2+1).
\end{equation}
In particular $f<g$. References to the \emph{smaller eigenspace} below always mean the smaller of the two nonconstant eigenspaces, namely the $r$-eigenspace. The constant eigenspace has dimension one, but projecting onto it would identify all vertices.

\section{The projection matrix and the independence number}
\subsection{The orthogonal projection}
Let $V_r$ be the $r$-eigenspace. The orthogonal projection onto it is
\begin{equation}\label{eq:E}
\boxed{E=\frac{2}{q(q+1)}A+\frac1qI-\frac1{q^2}J.}
\end{equation}
To check the formula, apply it to the three perpendicular eigenspaces. On constants its multiplier is
\[
\frac{2k}{q(q+1)}+\frac1q-q
=\frac{q^2-1}{q}+\frac1q-q=0.
\]
On $V_r$ the matrix $J$ vanishes and the multiplier is
\[
\frac{2r}{q(q+1)}+\frac1q=\frac{q-1}{q}+\frac1q=1.
\]
On the $s$-eigenspace it is $2s/(q(q+1))+1/q=-1/q+1/q=0$. Thus $E$ keeps exactly the $V_r$ component of a vector and removes the other two components. Consequently
\begin{equation}\label{eq:projectionproperties}
E^T=E,\quad E^2=E,\quad E\one=0,\quad \rank E=f=q(q-1).
\end{equation}
The projection is orthogonal because the decomposition into the three eigenspaces is perpendicular.

Its entries are
\begin{equation}\label{eq:Eentries}
E_{xy}=\begin{cases}
\dfrac{q-1}{q^2},&x=y,\\[5pt]
\dfrac{q-1}{q^2(q+1)},&x\sim y,\\[5pt]
-\dfrac1{q^2},&x\ne y,\ x\not\sim y.
\end{cases}
\end{equation}
For example, the edge entry is
\[
\frac{2}{q(q+1)}-\frac1{q^2}
=\frac{2q-(q+1)}{q^2(q+1)}
=\frac{q-1}{q^2(q+1)}.
\]

\subsection{Why columns have these dot products}
Let $e_x\in\R^n$ be the standard coordinate vector with one in position $x$ and zeros elsewhere. Define
\begin{equation}\label{eq:vx}
v_x=Ee_x.
\end{equation}
The vector $v_x$ is exactly column $x$ of $E$. Every such vector belongs to $V_r$. A dot product of column vectors is $a\cdot b=a^Tb$, so
\begin{align}
v_x\cdot v_y
&=(Ee_x)^T(Ee_y)\notag\\
&=e_x^TE^TEe_y\notag\\
&=e_x^TE^2e_y&&\text{since }E^T=E\notag\\
&=e_x^TEe_y&&\text{since }E^2=E\notag\\
&=E_{xy}.\label{eq:columnproducts}
\end{align}
The last step simply selects an entry: $Ee_y$ selects column $y$, and multiplication on the left by $e_x^T$ selects coordinate $x$ of that column. Equivalently,
\[
\sum_z(v_x)_z(v_y)_z
=\sum_zE_{zx}E_{zy}
=\sum_zE_{xz}E_{zy}=(E^2)_{xy}=E_{xy}.
\]
Thus the entries in~\eqref{eq:Eentries} are actual Euclidean dot products of explicitly specified real vectors.

\subsection{Both bounds for the independence number}
Proposition~\ref{prop:degree} gives an independent set of size $q$, so $\alpha(H_q)\ge q$. For the reverse inequality, let $x_1,\ldots,x_m$ be an independent set. Equations~\eqref{eq:Eentries} and~\eqref{eq:columnproducts} give
\[
\|v_{x_i}\|^2=\frac{q-1}{q^2},\qquad
v_{x_i}\cdot v_{x_j}=-\frac1{q^2}\quad(i\ne j).
\]
Therefore
\begin{align*}
0\le\left\|\sum_{i=1}^m v_{x_i}\right\|^2
&=\sum_i\|v_{x_i}\|^2+2\sum_{i<j}v_{x_i}\cdot v_{x_j}\\
&=\frac{m(q-1)}{q^2}-\frac{m(m-1)}{q^2}
=\frac{m(q-m)}{q^2}.
\end{align*}
For $m>0$ this forces $m\le q$; the empty set is harmless. Hence
\begin{equation}\label{eq:alphaH}
\boxed{\alpha(H_q)=q.}
\end{equation}
The upper bound just obtained is the Hoffman bound specialized to this graph: $n(-s)/(k-s)=q$; see Hoffman~\cite{Hoffman1970}. Our squared-norm calculation proves it directly in the present setting.

The same independent-set partition gives $\chi(H_q)=q^2$: at least $q^3/q=q^2$ colors are needed, and coloring each $I_{a,b}$ with its own color attains this number.

\section{The explicit two-distance embedding}
\subsection{From \texorpdfstring{$n$}{n} coordinates to \texorpdfstring{$f$}{f} coordinates}
The vectors $v_x$ in~\eqref{eq:vx} have $n$ entries, but they all belong to the same $f$-dimensional subspace $V_r$. Choose an orthonormal basis $b_1,\ldots,b_f$ of this subspace. Each $b_i$ is itself a vector in $\R^n$, and
\[
b_i\cdot b_j=\begin{cases}1,&i=j,\\0,&i\ne j.\end{cases}
\]
Write $v_x=\sum_i\alpha_i b_i$. Taking its dot product with $b_j$ gives
\[
b_j\cdot v_x=\sum_i\alpha_i(b_j\cdot b_i)=\alpha_j.
\]
The expansion is therefore
\begin{equation}\label{eq:basisexpansion}
v_x=\sum_{i=1}^f(b_i\cdot v_x)b_i.
\end{equation}
Define its coordinate vector in this basis by
\begin{equation}\label{eq:ux}
\boxed{u_x=(b_1\cdot v_x,\ldots,b_f\cdot v_x)\in\R^f.}
\end{equation}
There is no additional scaling factor in this definition.

\subsection{Precisely why inner products are preserved}
Put $\alpha_i=b_i\cdot v_x$ and $\beta_i=b_i\cdot v_y$. Expanding~\eqref{eq:basisexpansion} for both vectors yields
\begin{align*}
v_x\cdot v_y
&=\left(\sum_i\alpha_i b_i\right)\cdot
  \left(\sum_j\beta_j b_j\right)\\
&=\sum_i\sum_j\alpha_i\beta_j(b_i\cdot b_j)
=\sum_i\alpha_i\beta_i.
\end{align*}
The terms with $i\ne j$ vanish by orthogonality, while each term with $i=j$ retains its coefficient because $\|b_i\|=1$. The last sum is precisely the usual dot product of the coordinate tuples in~\eqref{eq:ux}. Hence
\begin{equation}\label{eq:uxgram}
\boxed{u_x\cdot u_y=v_x\cdot v_y=E_{xy}.}
\end{equation}
This equality is the exact meaning of preservation of inner products under passage to an orthonormal coordinate system.

There is also a useful entrywise description. Since $Eb_i=b_i$ and $E$ is symmetric,
\[
b_i\cdot v_x=b_i^TEe_x=(Eb_i)^Te_x=b_i^Te_x=b_i(x),
\]
where $b_i(x)$ means coordinate $x$ of the $n$-entry vector $b_i$. Therefore $u_x=(b_1(x),\ldots,b_f(x))$. The number $b_i(x)$ is at the same time the coefficient of $b_i$ in the expansion of $v_x$.

\subsection{A fully specified coordinate algorithm}
Equation~\eqref{eq:Eentries} already specifies every coordinate of $v_x$ directly from the graph. To obtain exactly $f$ coordinates per point without making an unspecified basis choice, choose an ordering of the field elements and then of $F^3$. Process the columns $v_x$ in that order by the following rule. If $b_1,\ldots,b_j$ have already been selected, form
\[
w=v_x-\sum_{i=1}^j(b_i\cdot v_x)b_i.
\]
If $w=0$, discard this column. Otherwise append $b_{j+1}=w/\|w\|$. The columns span $\im E$, so exactly $f$ basis vectors are selected. Then use~\eqref{eq:ux}. This is an exact construction using rational arithmetic and positive square roots; it requires no numerical eigensolver.

\subsection{Norms, distances, and dimension}
By~\eqref{eq:Eentries} and~\eqref{eq:uxgram},
\begin{equation}\label{eq:Hproducts}
\|u_x\|^2=\frac{q-1}{q^2},\qquad
u_x\cdot u_y=\begin{cases}
\dfrac{q-1}{q^2(q+1)},&x\sim y,\\[4pt]
-\dfrac1{q^2},&x\ne y,\ x\not\sim y.
\end{cases}
\end{equation}
Expanding coordinate squares gives $\|u-v\|^2=\|u\|^2+\|v\|^2-2u\cdot v$. Therefore
\begin{equation}\label{eq:Hdistances}
\boxed{\|u_x-u_y\|^2=\begin{cases}
\dfrac{2(q-1)}{q(q+1)},&x\sim y,\\[5pt]
\dfrac2q,&x\ne y,\ x\not\sim y.
\end{cases}}
\end{equation}
For example, the edge value is
\[
\frac{2(q-1)}{q^2}-\frac{2(q-1)}{q^2(q+1)}
=\frac{2(q-1)}{q(q+1)}.
\]
Both distances are positive, and the edge distance is strictly smaller. Thus all vertices are represented by different points. This is a two-distance embedding on the sphere of radius $\sqrt{q-1}/q$ in $\R^f$.

The columns of $E$ span $V_r$, so the $u_x$ span $\R^f$. Moreover,
\[
\sum_xv_x=E\sum_xe_x=E\one=0,
\]
and applying each coordinate functional gives $\sum_xu_x=0$. Thus the origin is their average. Their affine span contains the origin and hence equals their linear span, which has dimension $f$.

\section{The Taylor extension: definition and common-neighbor counts}
\subsection{Vertices, edges, and antipodal pairs}
For every vertex $x$ of $H_q$, introduce $x^+$ and $x^-$. Add $\infty^+$ and $\infty^-$. Define $T_q$ on these $2(n+1)$ vertices by the following rules, always for distinct vertices.
\begin{center}
\renewcommand{\arraystretch}{1.18}
\begin{tabular}{@{}ll@{}}
\toprule
Pair&Adjacency condition\\
\midrule
$x^+,y^+$ or $x^-,y^-$, $x\ne y$&$x\sim y$ in $H_q$\\
$x^+,y^-$, $x\ne y$&$x\not\sim y$ in $H_q$\\
$\infty^+,x^+$ or $\infty^-,x^-$&Always\\
$\infty^+,x^-$ or $\infty^-,x^+$&Never\\
$x^+,x^-$&Never\\
$\infty^+,\infty^-$&Never\\
\bottomrule
\end{tabular}
\end{center}
The pairs $\{x^+,x^-\}$ and $\{\infty^+,\infty^-\}$ are called antipodal pairs. Between any two different pairs there is a matching of two edges: original adjacency gives same-sign edges, original nonadjacency gives opposite-sign edges, and a pair involving the poles has same-sign edges. Thus every vertex has one neighbor in each of the other $n$ pairs. It follows that $T_q$ is $n$-regular and is a double cover of $K_{n+1}$ after the pairs are collapsed.

We call this graph the \emph{Taylor extension} of $H_q$. Equivalently, first adjoin a universal vertex to $H_q$, and then replace each vertex by a signed pair, using same-sign edges for adjacency and opposite-sign edges for distinct nonadjacency. Doubling $H_q$ alone would give $2n$ vertices; the poles are part of our definition.

If $C=J-I-A$ is the adjacency matrix of the complement of $H_q$, then in the order $\infty^+$, all $x^+$, all $x^-$, $\infty^-$, the adjacency matrix is
\begin{equation}\label{eq:Tmatrix}
A_T=\begin{pmatrix}
0&\one^T&0&0\\
\one&A&C&0\\
0&C&A&\one\\
0&0&\one^T&0
\end{pmatrix}.
\end{equation}

\subsection{Why \texorpdfstring{$k=2\mu$}{k = 2 mu} is decisive}
The edge rules can be applied to any graph, but uniform common-neighbor counts require additional structure. For an adjacent pair $\infty^+,x^+$, the common neighbors are exactly the $k$ vertices $z^+$ with $z\sim x$ in the input graph.

Now take an original nonedge $xy$. The pair $x^+,y^-$ is an edge in the extension. In the plus layer, a common neighbor corresponds to a neighbor of $x$ which is not a neighbor of $y$, giving $k-\mu$ choices. The minus layer gives another $k-\mu$, and neither pole is common. Thus this edge has $2(k-\mu)$ common neighbors. Equality with the pole-edge count requires
\[
2(k-\mu)=k\quad\Longleftrightarrow\quad k=2\mu.
\]
Our graph satisfies this condition.

We will also use the standard parameter identity, whose counting proof we include:
\begin{equation}\label{eq:srgcount}
k(k-\lambda-1)=(n-k-1)\mu.
\end{equation}
Fix $x$ and count edges from its $k$ neighbors to its $n-k-1$ distinct nonneighbors. Each neighbor has $k-\lambda-1$ neighbors in the latter set, while each nonneighbor has $\mu$ neighbors in the former set. This proves~\eqref{eq:srgcount}. Since $k=2\mu>0$,
\begin{equation}\label{eq:Bidentities}
B=2(k-\lambda-1),\qquad n+1=3k-2\lambda.
\end{equation}

\subsection{Every adjacent pair has \texorpdfstring{$k$}{k} common neighbors}
The pole-edge and opposite-layer cases were just counted. For the remaining type, suppose $x\sim y$ in $H_q$. Common neighbors of $x^+,y^+$ consist of their pole, $\lambda$ plus-layer vertices, and minus-layer vertices whose originals are adjacent to neither $x$ nor $y$.

Among the $n-2$ vertices other than $x,y$, exactly $2(k-1)-\lambda$ are adjacent to at least one endpoint. Therefore the number in the minus layer is $n-2k+\lambda$. The total is
\[
1+\lambda+n-2k+\lambda=n+1-2k+2\lambda=k,
\]
where the last step uses~\eqref{eq:Bidentities}. Reversing all signs covers the equivalent cases. Thus every adjacent pair in $T_q$ has precisely $k$ common neighbors.

\subsection{Common neighbors of nonadjacent pairs}
Antipodal vertices have no common neighbor: every vertex outside their pair is adjacent to exactly one of them. All other nonadjacent pairs, up to reversing signs, are in the following table.
\begin{center}
\renewcommand{\arraystretch}{1.2}
\begin{tabular}{@{}lll@{}}
\toprule
Pair&Condition in $H_q$&Common neighbors\\
\midrule
$\infty^+,x^-$&Any $x$&$n-k-1=B$\\
$x^+,y^+$&$x\ne y$, $x\not\sim y$&$n-1-2k+2\mu=B$\\
$x^+,y^-$&$x\sim y$&$2(k-\lambda-1)=B$\\
\bottomrule
\end{tabular}
\end{center}
For the first row, a common neighbor must be $z^+$, where $z$ is a distinct nonneighbor of $x$ in $H_q$. For the second row, there are $\mu$ common plus-layer neighbors, one common pole, and $n-2-2k+\mu$ common minus-layer neighbors. The last count excludes the endpoints and the union of their original neighbor sets. For the third row, each layer contributes $k-\lambda-1$ common neighbors: the subtraction of one removes the other endpoint. The displayed simplifications use $k=2\mu$ and~\eqref{eq:Bidentities}.

\section{Distance-regularity and the intersection array}
Graph distance is the length of a shortest edge path. It is distinct from the Euclidean distances to be constructed later. A connected graph is \emph{distance-regular} if, for every pair $w,z$ at distance $i$, the numbers of neighbors of $z$ at distances $i-1,i,i+1$ from $w$ depend only on $i$. These three numbers are denoted $c_i,a_i,b_i$.

Since $B>0$, every nonadjacent, nonantipodal pair in $T_q$ has graph distance two. An antipodal pair has no common neighbor and hence has distance greater than two. Its distance is exactly three: choose a neighbor $z$ of one member $w$. The matching rule shows that $z$ is a nonantipodal nonneighbor of the other member $w'$. There is a path of length two from $z$ to $w'$, producing a path of length three from $w$ to $w'$.

Fix a base vertex $w$. The distance layers have sizes $1,n,n,1$: the last layer is its antipode $w'$, and each other antipodal pair contributes one neighbor and one distance-two vertex.

A vertex in layer one has one neighbor in layer zero and $k$ in layer one, since adjacent pairs have $k$ common neighbors. Its remaining $n-k-1=B$ neighbors lie in layer two.

A vertex in layer two has $B$ neighbors in layer one, since it has $B$ common neighbors with $w$. It is adjacent to $w'$: it is not adjacent to $w$ and must be adjacent to one member of that antipodal pair. Its remaining $n-B-1=k$ neighbors lie in layer two. Finally, $w'$ is adjacent to every vertex in layer two.

The full table is
\begin{center}
\renewcommand{\arraystretch}{1.2}
\begin{tabular}{@{}ccccc@{}}
\toprule
Distance $i$&Layer size&Closer: $c_i$&Same: $a_i$&Farther: $b_i$\\
\midrule
0&1&0&0&$n$\\
1&$n$&1&$k$&$B$\\
2&$n$&$B$&$k$&1\\
3&1&$n$&0&0\\
\bottomrule
\end{tabular}
\end{center}
These values do not depend on the choice of $w$ or of the vertex in a layer. This proves distance-regularity. The intersection array lists the outward counts followed by the inward counts:
\begin{equation}\label{eq:intersection}
\boxed{\{b_0,b_1,b_2;\ c_1,c_2,c_3\}
=\left\{q^3,\frac{(q-1)(q^2+1)}2,1;\ 1,\frac{(q-1)(q^2+1)}2,q^3\right\}.}
\end{equation}
A distance-regular graph with array $\{N,M,1;1,M,N\}$ is called a Taylor graph~\cite[Sections 1.2.7 and 8.10.4]{BVM2022}. Notice that the array parameter $B$ is not the original strongly regular parameter $\mu$.

\section{The three-distance embedding of the Taylor extension}
\subsection{Choosing the extra coordinate}
Use the vectors $u_x\in\R^f$ from~\eqref{eq:ux}, with dot products~\eqref{eq:Hproducts}. Introduce a positive height $h$ and set
\[
p_{x^+}=(u_x,h),\qquad p_{x^-}=(-u_x,-h).
\]
These points have $f+1$ coordinates. Within either layer, distances are unchanged because the last coordinates cancel. The squared edge length remains
\[
\ell^2=\frac{2(q-1)}{q(q+1)}.
\]
For distinct original nonneighbors $x,y$, the Taylor rules require $x^+,y^-$ to be adjacent. Their squared distance is
\begin{align*}
\|p_{x^+}-p_{y^-}\|^2
&=\|u_x+u_y\|^2+(2h)^2\\
&=\frac{q-1}{q^2}+\frac{q-1}{q^2}-\frac2{q^2}+4h^2\\
&=\frac{2(q-2)}{q^2}+4h^2.
\end{align*}
Equating this to $\ell^2$ gives
\[
4h^2=\frac{2(q-1)}{q(q+1)}-\frac{2(q-2)}{q^2}
=\frac4{q^2(q+1)}.
\]
Therefore the required positive height is
\begin{equation}\label{eq:height}
\boxed{h=\frac1{q\sqrt{q+1}}.}
\end{equation}

\subsection{Placing the poles}
Place the upper pole at $p_{\infty^+}=(0,\ldots,0,t)$ with $t>h$. It must have squared distance $\ell^2$ from every plus-layer point. Thus
\[
\frac{q-1}{q^2}+(t-h)^2=\frac{2(q-1)}{q(q+1)},
\]
which gives
\[
(t-h)^2=\frac{(q-1)^2}{q^2(q+1)}.
\]
Taking the positive root,
\[
t=h+\frac{q-1}{q\sqrt{q+1}}=\frac1{\sqrt{q+1}}.
\]
Place the lower pole at its negative. The complete construction is
\begin{equation}\label{eq:allcoordinates}
\boxed{\begin{aligned}
p_{x^+}&=\left(u_x,\frac1{q\sqrt{q+1}}\right),&
p_{x^-}&=\left(-u_x,-\frac1{q\sqrt{q+1}}\right),\\[4pt]
p_{\infty^+}&=\left(0,\ldots,0,\frac1{\sqrt{q+1}}\right),&
p_{\infty^-}&=\left(0,\ldots,0,-\frac1{\sqrt{q+1}}\right).
\end{aligned}}
\end{equation}
The first $f$ coordinates of each pole are zero. Every finite point has squared norm
\[
\frac{q-1}{q^2}+\frac1{q^2(q+1)}=\frac1{q+1},
\]
and the poles have the same squared norm. Hence all $2(q^3+1)$ points lie on a sphere of radius $1/\sqrt{q+1}$. Antipodal vertices are represented by opposite vectors.

\subsection{Verification of all dot products}
Put $\rho=1/(q(q+1))$. Same-layer dot products for distinct originals are
\[
u_x\cdot u_y+h^2=\begin{cases}
\dfrac{q-1}{q^2(q+1)}+\dfrac1{q^2(q+1)}=\rho,&x\sim y,\\[5pt]
-\dfrac1{q^2}+\dfrac1{q^2(q+1)}=-\rho,&x\not\sim y.
\end{cases}
\]
The computation is identical in the minus layer, because both full vectors have been negated. For opposite-layer points with distinct originals, negating one full vector reverses the dot product. Thus original nonedges give $+\rho$ and original edges give $-\rho$, exactly in accordance with Taylor adjacency.

For poles and finite points,
\[
p_{\infty^+}\cdot p_{x^+}
=\frac1{\sqrt{q+1}}\frac1{q\sqrt{q+1}}=\rho,
\qquad
p_{\infty^+}\cdot p_{x^-}=-\rho.
\]
Reversing both signs gives the lower-pole cases. For antipodal pairs, the dot product is the negative of the squared norm, namely $-1/(q+1)$. All distinct pairs are therefore covered by
\begin{equation}\label{eq:Tproducts}
\boxed{p_a\cdot p_b=\begin{cases}
\dfrac1{q(q+1)},&a\sim b,\\[5pt]
-\dfrac1{q(q+1)},&a,b\text{ nonadjacent and nonantipodal},\\[5pt]
-\dfrac1{q+1},&a,b\text{ antipodal}.
\end{cases}}
\end{equation}

\subsection{The three distances}
Since the squared norm is $1/(q+1)$, the squared distance is $2/(q+1)-2p_a\cdot p_b$. Thus
\begin{center}
\renewcommand{\arraystretch}{1.55}
\begin{tabular}{@{}ccc@{}}
\toprule
Graph distance&Relationship&Squared Euclidean distance\\
\midrule
1&Adjacent&$\dfrac{2(q-1)}{q(q+1)}$\\
2&Nonadjacent, nonantipodal&$\dfrac2q$\\
3&Antipodal&$\dfrac4{q+1}$\\
\bottomrule
\end{tabular}
\end{center}
In particular, the three distances themselves are
\begin{equation}\label{eq:Tdistances}
\boxed{\sqrt{\frac{2(q-1)}{q(q+1)}},\qquad
\sqrt{\frac2q},\qquad\frac2{\sqrt{q+1}}.}
\end{equation}
They are strictly increasing for $q\ge3$. Hence all points are different and the pairs at the smallest distance are exactly the edges of $T_q$. In general these three squared distances are not proportional to $1,2,3$; that additional feature occurs for $q=3$.

For completeness, the pole distances can also be checked directly:
\[
\|p_{\infty^+}-p_{x^+}\|^2
=\frac{q-1}{q^2}+\frac{(q-1)^2}{q^2(q+1)}
=\frac{2(q-1)}{q(q+1)},
\]
\[
\|p_{\infty^+}-p_{x^-}\|^2
=\frac{q-1}{q^2}+\frac{(q+1)^2}{q^2(q+1)}
=\frac2q.
\]
Each antipodal distance is twice the common radius, giving squared distance $4/(q+1)$.

\subsection{Exact dimension and optional unit-edge normalization}
The original points $u_x$ are centered and span $\R^f$. Their differences also span it: for every $x$,
\[
u_x=\frac1n\sum_y(u_x-u_y).
\]
Differences within the plus layer therefore give $f$ independent horizontal directions. The difference between a pole and a plus-layer point has nonzero last coordinate, since $t-h>0$. It supplies an independent direction. Thus the affine dimension is exactly $f+1=q(q-1)+1$.

Our coordinates have retained the unscaled normalization throughout. If unit edges are desired afterwards, multiply all coordinates by
\[
\sqrt{\frac{q(q+1)}{2(q-1)}}.
\]
The three distances then become $1$, $\sqrt{(q+1)/(q-1)}$, and $\sqrt{2q/(q-1)}$. No such rescaling is needed for any of the arguments in this paper.

\section{Independence and the chromatic lower bound for \texorpdfstring{$T_q$}{Tq}}
\subsection{An independent set of size \texorpdfstring{$q+1$}{q + 1}}
Let $I$ be one of the independent $q$-sets $I_{a,b}$ in $H_q$. Then
\[
\{\infty^-\}\cup\{x^+:x\in I\}
\]
is independent in $T_q$. The plus layer retains the original nonadjacencies, and the lower pole has no plus-layer neighbors. Therefore $\alpha(T_q)\ge q+1$.

\subsection{The matching upper bound}
If an independent set contains an antipodal pair, it contains no third vertex: every vertex outside the pair is adjacent to exactly one member. Such an independent set has size at most two.

Otherwise, every distinct pair in an independent set is nonadjacent and nonantipodal. For a set $\{a_1,\ldots,a_m\}$ of this kind,~\eqref{eq:Tproducts} gives
\[
\|p_{a_i}\|^2=\frac1{q+1},\qquad
p_{a_i}\cdot p_{a_j}=-\frac1{q(q+1)}\quad(i\ne j).
\]
Expanding the squared norm of the sum,
\begin{align*}
0\le\left\|\sum_{i=1}^m p_{a_i}\right\|^2
&=\frac{m}{q+1}-\frac{m(m-1)}{q(q+1)}\\
&=\frac{m(q+1-m)}{q(q+1)}.
\end{align*}
For $m>0$, this forces $m\le q+1$. Combining this with the antipodal-pair case and the constructed lower bound proves
\begin{equation}\label{eq:alphaT}
\boxed{\alpha(T_q)=q+1=\alpha(H_q)+1.}
\end{equation}

\subsection{A chromatic lower bound}
Every color class in a proper coloring is independent, so it contains at most $q+1$ vertices. Consequently
\begin{equation}\label{eq:chiT}
\boxed{\chi(T_q)\ge\frac{2(q^3+1)}{q+1}=2(q^2-q+1).}
\end{equation}
Only this lower bound is asserted here; no optimal coloring of the extension is needed.

\subsection{Proof of the ball-packing theorem}\label{sec:packingproof}
\begin{proof}[Proof of Theorem~\ref{thm:packing}]
Let $d=q(q-1)+1=q^2-q+1$. Use the points $p_a\in\R^{d}$ defined in~\eqref{eq:allcoordinates}. Their minimum pairwise distance, attained precisely on the edges of $T_q$, is
\[
\ell_q=\sqrt{\frac{2(q-1)}{q(q+1)}}.
\]
At every point $p_a$, place a closed ball with radius
\begin{equation}\label{eq:packingradius}
r_q=\frac{\ell_q}{2}=\sqrt{\frac{q-1}{2q(q+1)}}.
\end{equation}
For distinct $a,b$, the center distance is at least $2r_q$, so the interiors of the two balls are disjoint. Equality holds exactly when $a\sim b$ in $T_q$, while every nonedge has strictly larger center distance. Thus these $2(q^3+1)$ equal-radius balls form a packing whose contact graph is exactly $T_q$.

Every independent set in this contact graph has at most $q+1$ vertices by~\eqref{eq:alphaT}. A proper coloring with $c$ colors therefore satisfies
\[
2(q^3+1)\le c(q+1).
\]
Since $q^3+1=(q+1)(q^2-q+1)$, we obtain
\[
c\ge 2(q^2-q+1)=2d.
\]
This particular congruent-ball packing proves $\chicong(d)\ge2d$. Allowing arbitrary radii can only enlarge the class of packings, giving $\chiball(d)\ge\chicong(d)$. These are the inequalities in~\eqref{eq:packing-intro}.

If unit-radius balls are preferred, dilate the entire packing by $1/r_q$. The new centers are $p_a/r_q$, the balls have radius one, and the contact graph is unchanged. The exact affine dimension of the centers is $d$, as proved in the preceding section.
\end{proof}

\section{Examples and the case \texorpdfstring{$q=3$}{q = 3}}
\subsection{Recovering the 27-vertex adjacency rule}
For $q=3$, take $F=\F_3$ and $\delta{}=2=-1$. Then
\[
D=(a-A)^2+(b-B)^2.
\]
When $(a,b)\ne(A,B)$, Lemma~\ref{lem:anisotropic} gives $D\ne0$, and every nonzero element of $\F_3$ has square one. Thus $D^2=1$, while $4\delta{}=2$. The edge-test expression becomes $2S^2-1$. It equals $2$, a nonsquare, when $S=0$, and equals $1$, a square, when $S\ne0$. Together with the same-horizontal-coordinate case, this gives
\[
(a,b,c)\sim(A,B,C)
\quad\Longleftrightarrow\quad
(a,b)\ne(A,B)\text{ and }c-C+bA-aB\ne0.
\]
This simplification is special to $q=3$; the general construction uses the square test in~\eqref{eq:adjacency}.

\subsection{All spectral and geometric constants for \texorpdfstring{$q=3$}{q = 3}}
The graph $H_3$ has parameters $(27,16,10,8)$, spectrum $16^1,4^6,(-2)^{20}$, and projection
\[
E=\frac16A+\frac13I-\frac19J.
\]
The column vectors, and their six-coordinate versions $u_x$, have dot products
\[
\|u_x\|^2=\frac29,\qquad
u_x\cdot u_y=\begin{cases}1/18,&x\sim y,\\-1/9,&x\ne y,\ x\not\sim y.\end{cases}
\]
The distances are $1/\sqrt3$ on edges and $\sqrt{2/3}$ on distinct nonedges. Its independence number is three.

The Taylor extension has 56 vertices, degree 27, and intersection array
\[
\{27,10,1;\ 1,10,27\}.
\]
It is the Gosset graph; the identification with the Schl\"afli local graph and this intersection array is standard~\cite{BVM2022}. Its embedding in $\R^7$ is
\[
p_{x^+}=(u_x,1/6),\qquad p_{x^-}=(-u_x,-1/6),
\]
with poles $(0,\ldots,0,\pm1/2)$. Every point has squared norm $1/4$. The three distances are
\[
\frac1{\sqrt3},\qquad\sqrt{\frac23},\qquad1.
\]
The independent-set argument gives $\alpha(T_3)=4$ and hence $\chi(T_3)\ge56/4=14$. In packing notation, Theorem~\ref{thm:packing} gives $\chiball(7)\ge\chicong(7)\ge14$, witnessed by 56 congruent balls. At the unscaled centers above, their radius is $r_3=1/(2\sqrt3)$. Multiplication of all coordinates by $\sqrt3$ would give distances $1,\sqrt2,\sqrt3$, but the unscaled coordinates are the ones used throughout our proofs.

\subsection{Further small prime powers}
The first four parameter sets are listed below. The dimension is the exact affine dimension of the embedding constructed here.
\begin{center}
\renewcommand{\arraystretch}{1.15}
\begin{tabular}{@{}rrrrrrr@{}}
\toprule
$q$&$|H_q|$&$k$&$\lambda$&$\mu$&$\dim H_q$&$\alpha(H_q)$\\
\midrule
3&27&16&10&8&6&3\\
5&125&72&45&36&20&5\\
7&343&192&116&96&42&7\\
9&729&400&235&200&72&9\\
\bottomrule
\end{tabular}
\end{center}
\begin{center}
\renewcommand{\arraystretch}{1.15}
\begin{tabular}{@{}rrrrrr@{}}
\toprule
$q$&$|T_q|$&Degree&$B$&$\dim T_q$&$\alpha(T_q)$\\
\midrule
3&56&27&10&7&4\\
5&252&125&52&21&6\\
7&688&343&150&43&8\\
9&1460&729&328&73&10\\
\bottomrule
\end{tabular}
\end{center}
For $q=5$, one may take $\delta{}=2$ in $\F_5$. A vertex $(A,B,C)$ is adjacent to $(0,0,0)$ precisely when
\[
3C^2-(A^2-2B^2)^2\in\{1,4\}.
\]
For instance, when $(A,B)=(1,0)$, the allowed values are $C=0,2,3$, illustrating the $(q+1)/2=3$ neighbors above each different horizontal pair.

For $q=9$, take $F=\F_3[\omega]$ with $\omega^2=2$. Its nonzero squares are $1,2,\omega,2\omega$, so $\delta{}=1+\omega$ is a valid nonsquare. The same triple rule constructs all 729 vertices and their edges using this field arithmetic. The extra symbol $\theta$ used in the strong-regularity proof then generates a quadratic extension of this nine-element field; it is separate from the symbol $\omega$ used to construct $F$ itself.

\end{document}